\documentclass[10pt,leqno]{amsart}
\usepackage{indentfirst, amssymb,amsmath,amsthm}
\newtheorem*{theoA}{Theorem A}
\newtheorem*{theoB}{Theorem B}
\newtheorem*{theoC}{Theorem C}
\newtheorem*{theoD}{Theorem D}

\newtheorem{theo}{Theorem}[section]
\newtheorem{lem}{Lemma}[section]
\newtheorem{cor}{Corollary}[section]

\newtheorem{exm}{Example}[section]
\newtheorem{defi}{Definition}[section]
\newtheorem{rem}{Remark}[section]
\newcommand{\ol}{\overline}
\newcommand{\be}{\begin{equation}}
\newcommand{\ee}{\end{equation}}
\newcommand{\beas}{\begin{eqnarray*}}
	\newcommand{\eeas}{\end{eqnarray*}}
\newcommand{\bea}{\begin{eqnarray}}
\newcommand{\eea}{\end{eqnarray}}
\newcommand{\lra}{\longrightarrow}
\numberwithin{equation}{section}
\begin{document}
	 { T\MakeLowercase  {\large his paper has been accepted for publication in \bf \MakeUppercase {Studia Universitatis Babes-Bolyai Mathematica} on 17.09.2020.}}
\title[ \large L\MakeLowercase {\large inear Delay-Differential operator of a meromorphic function...}]{ L\MakeLowercase  {\large inear Delay-Differential operator of a meromorphic function sharing two sets or small function together with values with its $c$-shift or $q$-shift}}
\numberwithin {equation}{section}
\date{}
\author{Arpita Roy\;\; and \;\;Abhijit Banerjee}\\

\date{}
\address{ Department of Mathematics, University of Kalyani, West Bengal 741235, India.}
\email{arpita140793@gmail.com}
\address{Department of Mathematics, University of Kalyani, West Bengal 741235, India.}
\email{abanerjee\_kal@yahoo.co.in, abanerjeekal@gmail.com}
\
\maketitle
\let\thefootnote\relax
\footnotetext{Type set by \AmS -\LaTeX}
\begin{abstract}
	 The paper is devoted to study the uniqueness problem of linear delay-differential operator of a meromorphic function sharing two sets or small function together with values with its $c$-shift and $q$-shift operator. Results of this paper drastically improve two recent results of Meng-Liu [J. Appl. Math. Inform. 37(1-2)(2019), 133-148] and Qi-Li-Yang [Comput. Methods Funct. Theo., 18(2018), 567-582]. In addition to this, one of our result  improve and extend that of Qi-Yang [Comput. Methods Funct. Theo., 20(2020), 159-178].
 \end{abstract}
\section{Introduction Definitions and Results}
 Throughout the paper we use standard notations of Nevanlinna theory as stated in \cite{Hayman} and by any meromorphic function $f$ we always mean that it is defined on $\mathbb{C}$.  Let $f$ and $g$ be such two non constant meromorphic functions. 
For $a\in \mathbb{C}\cup \{\infty\}$, the following two quantities \[\delta(a;f)= 1- \limsup\limits_{r\lra \infty}\frac{N(r,a;f)}{T(r,f)}= \liminf\limits_{r\lra \infty}\frac{m(r,a;f)}{T(r,f)}\] and $$\Theta(a;f)= 1- \limsup\limits_{r\lra \infty}\frac{\ol N(r,a;f)}{T(r,f)}$$ are respectively known as Nevanlinna deficiency and ramification index of the value $a$.
\par In the beginning of the nineteenth century R. Nevanlinna inaugurated the value distribution theory with his famous Five value and Four value theorems which can be considered as the backbone of the modern uniqueness theory. Illuminated by these two basic results initially the research were performed on the value sharing of meromorphic functions.  After five decades, uniqueness theory moved to a new direction led by F. Gross \cite{Gross_LNM}, who transformed the traditional value sharing problem to a more general set up namely to the shared set problems. Now we recall the definition of set sharing.

\begin{defi}Let $S$ be a set of distinct elements of $\mathbb{C}\cup\{\infty\}$.
	For a  non-constant meromorphic function $f$, let $E_{f}(S)=\bigcup_{a\in S}\{(z,p)\in\mathbb{C}\times\mathbb{N}: f(z)=a\; with\; multiplicity\; p\}$ $\left(\ol  E_{f}(S)=\bigcup_{a\in S}\{(z,1)\in\mathbb{C}\times\mathbb{N}: f(z)=a\}\right) $. Then we say  $f$, $g$ share the set $S$ CM(IM) if $E_{f}(S)=E_{g}(S)$ $\left(\ol  E_{f}(S)=\ol E_{g}(S)\right) $. 
\end{defi}
Evidently, if $S$ is a singleton, then it coincides with the traditional definition of CM(IM)  sharing of values, which are known to the readers. \par
In 2001, due to a revolutionary approach by Lahiri \cite{Lahiri_NMJ,Lahiri_CVTA}, the notion of weighted sharing of values or sets appeared in the literature  and expedite the research work there in. Though now-a-days the definition is widely circulated, we invoke the definition.
\begin{defi}\cite{Lahiri_NMJ,Lahiri_CVTA} Let $k$ be a non-negative integer or infinity. For $a\in\mathbb{C}\cup\{\infty\}$ we denote by $E_{k}(a;f)$ the set of all $a$-points of $f$, where an $a$-point of multiplicity $m$ is counted $m$ times if $m\leq k$ and $k+1$ times if $m>k$. If $E_{k}(a;f)=E_{k}(a;g)$, we say that $f,g$ share the value $a$ with weight $k$ and denote it by $(a,k)$. The IM and CM sharing corresponds to $(a,0)$ and $(a,\infty)$ respectively.\end{defi}

\begin{defi} \cite{Lahiri_NMJ} Let $S$ be a set of distinct elements of $\mathbb{C}\cup\{\infty\}$ and $k$ be a non-negative integer or $\infty$. We denote by $E_{f}(S,k)$ the set $\cup_{a\in{S}}E_{k}(a;f)$. Clearly $E_{f}(S)=E_{f}(S,\infty)$ and $\ol E_{f}(S)=E_{f}(S,0)$.\par If $E_{f}(S,k)=E_{g}(S,k)$, then we say that $f$, $g$ share the set $S$ with weight $k$ and write it as $f$, $g$ share $(S,k)$. \end{defi}	By $N(r,a;f\mid <m)$ we mean the counting function of those $a$-points of $f$ whose multiplicities are less than $m$ where each $a$-point is counted according to its multiplicity and by
	$\ol N(r,a;f\mid\geq m)$ we mean the counting function of those $a$-points of $f$ whose multiplicities are not less than $m$ where each $a$-point is counted ignoring multiplicity.
	We also denote by $N_{2}(r,a;f)$ the sum $\ol N(r,a;f)+\ol N(r,a;f\mid\geq 2)$. 

Usually, $S(r,f)$ denotes any quantity satisfying $S(r,f)= o(T(r,f))$ for all $r$ outside of a possible exceptional set of finite linear measure. Also $S_1(r,f)$ denotes any quantity satisfying $S_1(r,f)=o(T(r,f))$ for all $r$ on a set of logarithmic density 1, where the logarithmic density of a set $F$ is defined by
$$
\limsup_{r\to\infty}\frac{1}{\log r}\int_{[1,r]\cap
F}\frac{dt}{t}.
$$ \par

Throughout the paper for a positive integer $n$, $S_{1}$, $S^{*}_{1}$ and $S_{2}$ represents respectively the sets $\{1,\omega,\ldots,\omega^{n-1}\}$, $\{\alpha_{1}, \alpha_{2},\ldots, \alpha_{n}\}$ and $\{\infty\}$ , where $\omega =\cos\frac{2\pi}{n}+i\sin\frac{2\pi}{n}$ and $\alpha_{i}$, $i=1,2,\ldots,n$ are non zero constants. \par
Let $a_{t-1}(\neq 0)$, $a_{t-2}$, \ldots, $a_{0}$ and $C(\neq 0)$ be complex numbers. We define \bea\label{e1.1} P(z)=Cz Q(z)=Cz(a_{t-1}z^{t-1}+a_{t-2}z^{t-2}+ \ldots+a_{1}z+a_{0}). \eea 
 For the polynomial $P(z)$ as given in (\ref{e1.1}), let us define the following two functions:
\[
\chi^{t-1}_{_0}=
\begin{cases}
1, & if\; a_{0} \neq 0 \\
0, & if\; a_{0}=0
\end{cases} 
\] and 
\[
\mu^{t-1}_{_0}=
\begin{cases}
1, & if\; a_{0}=0, a_{1} \neq 0 \\
0, & otherwise.
\end{cases} 
\] 
In view of (\ref {e1.1}), corresponding to the set $S^{*}_{1}$, let us consider the polynomial $P_{*}(z)$ as follows:
\bea\label{e1.2}&& P_{*}(z)=Cz Q_{*}(z), \hspace{2mm}\text{where}\hspace{2mm}  C=\frac{1}{(-1)^{n+1}\alpha_{1}\alpha_{2} \ldots \alpha_{n}} \hspace{2mm}\text{and}\\&& Q_{*}(z)=\displaystyle\sum_{r=0}^{n-1}(-1)^{r}\sum \alpha_{1}\alpha_{2} \ldots \alpha_{r}\hspace{1mm} z^{n-r-1}\nonumber, \eea $\sum\alpha_{1}\alpha_{2} \ldots \alpha_{r}$ = sum of the products of the values $\alpha_{1},\alpha_{2}, \ldots ,\alpha_{n}$ taken $r$ into account. We also denote by $m_{1}$ and $m_{2}$ as the number of simple and multiple zeros of $Q_{*}(z)$ respectively.\par
Next we define linear shift operator, shift-differential operator and differential operator respectively as follows: 
 \beas\hspace{-1.0cc} L_{1}(f(z))=a_{k}f(z+c_{k})+a_{k-1}f(z+c_{k-1})+\ldots+a_{1}f(z+c_{1})+a_{0}f(z),\eeas \beas \hspace{2.2cc}L_{2}(f(z))=b_{s}f^{(s)}(z+c_{s})+b_{s-1}f^{(s-1)}(z+c_{s-1})+\ldots+b_{1}f^{\prime}(z+c_{1}), \eeas  \beas \hspace{-3.5cc} L_{3}(f(z))=d_{t}f^{(t)}(z)+d_{t-1}f^{(t-1)}(z)+\ldots+d_{1}f^{\prime}(z), \eeas where $a_{k}$, $b_{s}$ and $d_{t}$ are non zero and $k$, $s$, $t$ are natural numbers and all $c_{i}'s$ are non zero. For the sake of convenience we shall call $L_{2}(f(z))+L_{3}(f(z))$ as delay-differential operator which is denoted by $\Tilde{L}(f(z))$.  \par
As far as the knowledge of the authors are concerned, Qi-Li-Yang \cite{Qi-Li and Yang_CMFT} were the first authors who initiated two shared set problems for the derivative of a meromorphic function $f(z)$ with its shift $f(z+c)$ as follows:  
\begin{theoA}\cite{Qi-Li and Yang_CMFT}
Let $f(z)$ be a non constant meromorphic function of finite order, $n\geq 9$ be an integer and $a$ be a non zero complex constant. If $[f^{\prime}(z)]^{n}$ and $f^{n}(z+c)$ share $(a, \infty)$ and $(\infty, \infty)$, then $f^{\prime}(z)=tf(z+c)$, for a constant $t$ that satisfies $t^{n}=1$.
\end{theoA}
Recently employing the notion of weighted sharing, Meng-Liu \cite{Meng and Liu_JAMI} further investigated {\it Theorem A} to obtain the following result.
\begin{theoB}\cite{Meng and Liu_JAMI} 
Let $f(z)$ be a non constant meromorphic function of finite order, $n\geq10$ be an integer. If $[f^{\prime}(z)]^{n}$and $f^{n}(z+c)$ share $(1, 2)$ and $(\infty, 0)$, then $f^{\prime}(z)=tf(z+c)$, for a constant $t$ that satisfies $t^{n}=1$.
\end{theoB}
Considering $f(z)=e^{z}$ and $\omega =e^{-c}$ satisfying $\omega^{n}=1$, it is easy to see that $f^{\prime}$ and $f(z+c)$ share the sets $(S_{1}, \infty)$, $(\infty, \infty)$ and $f^{\prime}(z)=\omega f(z+c)$ for each $n$. So it is natural to conjecture that in {\it Theorem A} and {\it Theorem B} the cardinality of $n$ could further be reduced. To this end, we have performed our investigations and have been able to reduce the cardinality of $n$ in {\it Theorem B} up to $6$. In fact, we have proved our theorem for a more general setting $S_{1}^{*}$  rather tan to consider only the set $S_{1}$. 
\begin{theo}\label{t1}
Let $f(z)$ be a non constant meromorphic function of finite order such that $\Tilde{L}(f(z))$ and $f(z+c)$ share $(S_{1}^{*},  2)$ and $(S_{2}, 0)$. If \beas n > 2(\chi^{n-1}_{0}+\mu^{n-1}_{0}+m_{1}+2m_{2})+\frac{15}{(2n-3)} (\chi^{n-1}_{0}+m_{1}+m_{2}), \hspace{2mm} \text{then}\eeas 
\beas \prod_{i=1}^{n} (\Tilde{L}(f(z))-\alpha_{i}) \equiv \prod_{i=1}^{n} (f(z+c)-\alpha_{i}).\eeas

\end{theo}
\begin{rem}\label{r1.1}
	From the definitions, we easily can calculate the value of $\chi^{n-1}_{0}$, $\mu^{n-1}_{0}$, $m_{1}$ and $m_{2}$ for a particular set  $S_{1}^{*}$. Clearly for the set $S_{1}$, $\chi^{n-1}_{0}=0$; $\mu^{n-1}_{0}=0$; $m_{1}=0$ and $m_{2}=1$. Therefore in above theorem for the set $S_{1}$ if $n>4+ \frac{15}{(2n-3)}$ i.e., if $n\geq 6$ then $\Tilde{L}(f(z))=tf(z+c)$, for a constant $t$ that satisfies $t^{n}=1$. For a particular choices of coefficients of $\Tilde{L}(f(z))$ we can easily make $\Tilde{L}(f(z))=f^{\prime}$.
\end{rem}
Corresponding to $q$-shift Meng-Liu \cite{Meng and Liu_JAMI} also investigated the same result like {\it Theorem B} as follows :
\begin{theoC}\cite{Meng and Liu_JAMI} 
	Let $f(z)$ be a non constant meromorphic function of zero order, $n\geq10$ be an integer. If $[f^{\prime}(z)]^{n}$and $f^{n}(qz)$ share $(1, 2)$ and $(\infty, 0)$, then $f^{\prime}(z)=tf(qz)$, for a constant $t$ that satisfies $t^{n}=1$.
\end{theoC}
In connection to {\it Theorem C} below we present our result which improve the same.
\begin{theo}\label{t2}
	Let $f(z)$ be a non constant meromorphic function of zero order such that $\Tilde{L}(f(z))$ and $f(qz)$ share $(S_{1}^{*},  2)$ and $(S_{2}, 0)$.  If \beas n > 2(\chi^{n-1}_{0}+\mu^{n-1}_{0}+m_{1}+2m_{2})+\frac{15}{(2n-3)} (\chi^{n-1}_{0}+m_{1}+m_{2}) \hspace{2mm} \text{then}\eeas 
	\beas \prod_{i=1}^{n} (\Tilde{L}(f(z))-\alpha_{i}) \equiv \prod_{i=1}^{n} (f(qz)-\alpha_{i}).\eeas
\end{theo}

 In the next theorem we shall show that the lower bound of $n$ can further be reduced at the expense of allowing both the range sets $S_{1}^{*}$, $S_{2}$ to be shared CM. 
\begin{theo}\label{t3}
Let $f(z)$ be a non constant meromorphic function of finite order such that $\Tilde{L}(f(z))$ and $f(z+c)$ share $(S^{*}_{1}, \infty)$ and $(S_{2}, \infty)$ with $T(r,f)=N\left(r, \frac{1}{\Tilde{L}(f(z))}\right)+S(r,f)$ then for $n > 2(\chi^{n-1}_{0}+m_{1}+m_{2})+1$, \beas \prod_{i=1}^{n} (\Tilde{L}(f(z))-\alpha_{i}) \equiv \prod_{i=1}^{n} (f(z+c)-\alpha_{i}).\eeas
\end{theo}
\begin{rem}\label{r1.2} In connection of {\it Remark \ref{r1.1}}, for the set $S_{1}$ in {\it Theorem \ref{t3}} the result holds for $n\geq 4$.
\end{rem}
Our next theorem is analogous theorem of {\it Theorem \ref{t3}} corresponding to $q$-shift.
\begin{theo}\label{t4}
 Let $f(z)$ be a non constant meromorphic function of zero order such that $\Tilde{L}(f(z))$ and $f(qz)$ share $(S^{*}_{1}, \infty)$ and $(S_{2}, \infty)$.  If $ n > 2(\chi^{n-1}_{0}+m_{1}+m_{2})+1$ then 
\beas \prod_{i=1}^{n} (\Tilde{L}(f(z))-\alpha_{i}) \equiv \prod_{i=1}^{n} (f(qz)-\alpha_{i}).\eeas
\end{theo}
Recently, corresponding to {\it Theorem A}, Qi-Yang \cite{Qi and Yang_CMFT} obtained the value sharing problem for entire function as follows:
\begin{theoD}\cite{Qi and Yang_CMFT}
	Let $f(z)$ be a transcendental entire function of finite order and let $(a\neq0)\in\mathbb{C}$. If $f^{\prime}(z)$ and $f(z+c)$ share $(0, \infty)$ and $(a,0)$, then $f^{\prime}(z)\equiv f(z+c)$.
\end{theoD}
In view of {\it Theorem 1.1}, \cite{Qi and Yang_CMFT} we know that $f(z)$ actually becomes a transcendental entire function. Since we are dealing with $\Tilde{L}(f(z))$ instead of $f^{\prime}$, it will be reasonable to consider the above theorem for meromorphic function under small function sharing category. In this respect we prove the following theorem.
\begin{theo}\label{t5}
Let $f(z)$ be a transcendental meromorphic function of finite order and let $a(z)(\not\equiv0)\in S(f)$ be an entire function. If $\Tilde{L}(f(z))$ and $f(z+c)$ share $(0, \infty)$, $(\infty, \infty)$ and $(a(z),0)$ with $\Theta(0;f)+\Theta(\infty; f)> 0$, then $\Tilde{L}(f(z))\equiv f(z+c)$.
\end{theo}
From {\it Theorem \ref{t5}} we can immediately deduce the following corollary.
\begin{cor}
Let $f(z)$ be a transcendental entire function of finite order and let $a(z)(\not\equiv0)\in S(f)$. If $\Tilde{L}(f(z))$ and $f(z+c)$ share $(0, \infty)$ and $(a(z),0)$, then $\Tilde{L}(f(z))\equiv f(z+c)$.
\end{cor}
 Following example shows that in {\it Theorem \ref{t5}} the CM pole sharing can not be replaced by IM.
\begin{exm}\label{ex1}
Let $f(z)=\frac{2e^{2\sqrt{2}iz}-8e^{\sqrt{2}iz}+2}{(e^{\sqrt{2}iz}+1)^{2}}$ and $c=\sqrt{2}\pi$. Choose the coefficients of $\tilde{L}(f(z))$ in such a way that $\tilde{L}(f(z))=f^{\prime\prime}$. Then $\tilde{L}(f(z))\left(=\frac{24e^{\sqrt{2}iz}[e^{2\sqrt{2}iz}-4e^{\sqrt{2}iz}+1]}{(e^{\sqrt{2}iz}+1)^{4}}\right)$ and $f(z+c)$ share $(0, \infty)$, $(1, 0)$ and $(\infty, 0)$ and $\Theta(0;f)+\Theta(\infty; f)=\frac{1}{2}> 0$ but $\Tilde{L}(f(z))\not\equiv f(z+c)$.
\end{exm}
From the next example we can show that in Theorem \ref{t5} sharing of $0$ can not be replaced by sharing of a non zero value. 
\begin{exm}\label{ex2}
	Let $f(z)=(e^{\lambda z}-1)^{2}+1$. Choose $e^{\lambda c}=1$, $\displaystyle \sum_{i=1}^{s}b_{i}(2\lambda)^{i}e^{2\lambda c_{i}}+\displaystyle \sum_{i=1}^{t}d_{i}(2\lambda)^{i}=0$ and \; $\displaystyle \sum_{i=1}^{s}b_{i}(\lambda)^{i}e^{\lambda c_{i}}+\displaystyle \sum_{i=1}^{t}d_{i}(\lambda)^{i}=-\frac{1}{2}$. Then  $f(z+c)=(e^{\lambda z}-1)^{2}+1$ and $\Tilde{L}(f(z))=e^{\lambda z}$. Clearly $f(z+c)$ and $\Tilde{L}(f(z))$ share $(2, \infty)$, $(\infty, \infty)$ and $(1,0)$ with $\Theta(0;f)+\Theta(\infty; f)> 0$. But $\Tilde{L}(f(z))\neq f(z+c)$.	
\end{exm}
In {\it Theorem \ref{t5}}, sharing of the value $0$ can be removed at the cost of slightly manipulating the deficiency condition. In this respect, we state the following theorem for transcendental meromorphic function.
\begin{theo}\label{t6}
Let $f(z)$ be a transcendental meromorphic function of finite order and let $a(z)(\not\equiv0)\in S(f)$ be an entire function. If $\Tilde{L}(f(z))$ and $f(z+c)$ share  $(a(z),\infty)$ and $(\infty,\infty)$ with $\delta(0;f)> 0$, then $\Tilde{L}(f(z))\equiv f(z+c)$.	
\end{theo}

By an example we now show that $a(z)$ CM sharing can not be replaced by IM in {\it Theorem \ref{t6}}.
\begin{exm}\label{ex3}
	Let $f(z)=\frac{-2e^{z}-1}{e^{2z}}$ and $c=\pi i$. Choose $\Tilde{L}(f(z))=L_{3}(f(z))$ with\;$2\displaystyle\sum_{i=1}^{t}(-1)^{i+1}d_{i}=1$ and \;$\displaystyle\sum_{i=1}^{t}(-2)^{i}d_{i}=0$. Then $\Tilde{L}(f(z))=\frac{1}{e^{z}}$ and $f(z+c)=\frac{2e^{z}-1}{e^{2z}}$ share $(1,0)$, $(\infty, \infty)$ and $\delta(0;f)=\frac{1}{2}> 0$. Clearly $\Tilde{L}(f(z))\neq f(z+c)$.
\end{exm}
Our next example shows that $a(z)\not\equiv 0$ in {\it Theorem \ref{t6}} can not be dropped  as well as $(a(z),0)$ sharing in {\it Theorem \ref{t5}} can not be removed.
\begin{exm}\label{ex4}
	Let $f(z)=e^{\frac{\pi iz}{c}}$. Choose $\Tilde{L}(f(z))=f^{\prime}$. Then clearly $f(z+c)$ and $\Tilde{L}(f(z))$ share $(0, \infty)$, $(\infty, \infty)$  and $\delta(0;f)> 0$. But $\Tilde{L}(f(z))\neq f(z+c)$.
\end{exm}

Following two examples show that $\delta(0;f)> 0$ in {\it Theorem \ref{t6}} can not be removed.
\begin{exm}\label{ex5}
In {\it Example \ref{ex2}} though $f(z+c)$ and $\Tilde{L}(f(z))$ share $(2, \infty)$, $(\infty, \infty)$ but $\delta(0;f)= 0$. Here $\Tilde{L}(f(z))\neq f(z+c)$.	
\end{exm}
\begin{exm}\label{ex6}
Let $f(z)=\frac{e^{z}+z}{2}$ and $a(z)=z$. Choose $\Tilde{L}(f(z))=L_{3}(f(z))$ with $d_{1}=2c$ and $\displaystyle\sum_{j=2}^{t} d_{j}=2(e^{c}-c)$. Then $f(z+c)$ $\left(=\frac{e^{c}e^{z}+z+c}{2}\right)$ and $\Tilde{L}(f(z))$$(=e^{c}e^{z}+c)$ share $(a(z), \infty)$ and $(\infty, \infty)$ but $\delta(0;f)= 0$. Clearly $\Tilde{L}(f(z))\neq f(z+c)$.
\end{exm}

\section{Lemmas}
In this section some lemmas will be presented which will be needed in the sequel.

\begin{lem} \cite{Chiang and Feng_RJ} \label{l2.1}
	Let $f(z)$ be a meromorphic function of finite order $\rho$ and let $c\in\mathbb{C}\setminus\{0\}$ be fixed. Then, for each $\varepsilon > 0$, we have \beas T(r,f(z+c)) =T(r,f(z)) + O(r^{\rho-1+\varepsilon}) + O(\log r).\eeas 
\end{lem}
\begin{lem} \cite{Halburd and Korhonen_JMAA} \label{l2.2}
	Let $f(z)$ be a meromorphic function of finite order and $c \in\mathbb{C}\setminus\{0\}$. Then  
	\beas m\left(r,\frac{f(z+c)}{f(z)}\right)+m\left(r,\frac{f(z)}{f(z+c)}\right)= S(r,f).\eeas
\end{lem}
\begin{lem}\label{l2.3}\cite{Halburd-Korhonen and Tohge_TAMS} Let $f$ be a non-constant meromorphic function of finite order and $c\in \mathbb{C}$. Then 
	\beas N \left(r,\frac{1}{f(z+c)}\right)\leq N \left(r,\frac{1}{f(z)}\right)+S(r,f),\eeas \beas N \left(r,f(z+c)\right) \leq N\left(r,f(z)\right)+S(r,f),\eeas \beas \ol N \left(r,\frac{1}{f(z+c)}\right)\leq \ol N \left(r,\frac{1}{f(z)}\right)+S(r,f)\eeas and \beas \ol N \left(r,f(z+c)\right) \leq\ol N\left(r,f(z)\right)+S(r,f).\eeas \end{lem}	
\begin{lem}\cite{Barnett-Halburd-Korhonen and Morgan_PRSES} \label{l2.4}
	Let $f(z)$ be a meromorphic function of zero order and $q \in\mathbb{C}\setminus\{0\}$. Then  
	\beas m\left(r,\frac{f(qz)}{f(z)}\right)= S_{1}(r,f).\eeas
\end{lem}
\begin{lem} \cite{Zhang and Korhonen_JMAA} \label{l2.5}
Let $f(z)$ be a non constant zero order meromorphic function and $ q \in \mathbb{C} \setminus \{0\}$, then \beas T(r,f(qz)) = (1+o(1))T(r,f(z)) \eeas and \beas N(r,f(qz)) = (1+o(1))N(r,f(z))  \eeas on a set of lower logarithmic measure $1$. \end{lem}
Using {\it Lemma \ref{l2.4}} and {\it Lemma \ref{l2.5}} and by the help of simple transformation one can easily prove the next lemma.
\begin{lem}\label{l2.6} 
Let $f(z)$ be a meromorphic function of zero order and $q \in\mathbb{C}\setminus\{0\}$. Then  
\beas m\left(r,\frac{f(z)}{f(qz)}\right)= S_{1}(r,f).\eeas
\end{lem}

\begin{lem}\cite{Yang and Yi_Dordrecht} \label{l2.7} 
	Let $f(z)$ be a non constant meromorphic function in the complex plane, and let $R(f)=\frac{P(f)}{Q(f)}$, where \beas P(f)=\sum_{k=0}^{p}a_{k}(z)f^{k} \hspace{0.5cc} \text{and}\hspace{0.5cc} Q(f)=\sum_{j=0}^{q}b_{j}(z)f^{j}\eeas are two mutually prime polynomials in $f$. If the coefficients $a_{k}(z)$ for $k=0,1, \ldots, p$ and $b_{j}(z)$ for $j=0,1, \ldots, q$ are small functions of $f$ with $a_{p}(z)\not\equiv 0$ and $b_{q}(z)\not\equiv 0$, then \beas T(r, P(f))=\max\{p, q\} T(r,f)+S(r,f).\eeas
\end{lem}
\begin{lem}\cite{Li and Wang_JMAA}\label{l2.8}
	Suppose that $h$ is a non constant meromorphic function satisfying
\beas	N(r, h) + N\left(r, \frac{1}{h}\right) = S(r, h).\eeas
	Let $f = a_{0}h^{p} + a_{1}h^{p-1} + \ldots + a_{p}$, and $g = b_{0}h^{q} + b_{1}h^{q-1} + \ldots+ b_{q}$ be polynomials in $h$ with coefficients $a_{0}$, $a_{1}, \ldots$, $a_{p}$; $b_{0}$, $b_{1},\ldots$, $b_{q}$ being small
	functions of $h$ and $a_{0}b_{0}a_{p}\not\equiv 0$. If $q\leq p$, then $m\left(r, \frac{g}{f}\right) = S(r, h)$.
\end{lem}
\begin{lem}\cite{Lahiri and Dewan_KMJ}\label{l2.9} If $N(r,0;f^{(k)}\mid f\not=0)$ denotes the counting function of those zeros of $f^{(k)}$ which are not the zeros of $f$, where a zero of $f^{(k)}$ is counted according to its multiplicity then $$N(r,0;f^{(k)}\mid f\not=0)\leq k\ol N(r,\infty;f)+N(r,0;f\mid <k)+k\ol N(r,0;f\mid\geq k)+S(r,f).$$\end{lem}
\begin{lem}\label{l2.10} Let $F$ be a meromorphic function. Then $$\ol N(r,1;F\mid\geq k+1)\leq\frac{1}{k}\{\ol N(r,0;F)+\ol N(r,\infty;F)\}+S(r,F).$$\end{lem}
Since the proof is straight forward, it is omitted.
\begin{lem}\label{l2.11}\cite{Banerjee_CMJ} Let $F$, $G$ be two meromorphic functions sharing $(1,2)$ and $(\infty,k)$, where $0\leq k\leq\infty$. Then one of the following cases holds\\
	\beas (i)\;T(r,F)+T(r,G)&\leq& 2\{N_{2}(r,0;F)+N_{2}(r,0;G)+\ol N(r,\infty;F)+\ol N(r,\infty;G)\;\;\;\;\;\;\;\;\;\;\;\;\;\;\;\;\;\;\;\;\;\;\;\;\;\\& & +\ol N_{*}(r,\infty;F,G)\}+S(r,F)+S(r,G),\eeas where $\ol N_{*}(r,\infty;f,g)$ is the reduced counting function of those poles of $F$ whose multiplicities differ from the multiplicities of the corresponding poles of $G$,\\(ii)$\;F\equiv G$,\\(iii)$\;FG\equiv 1$.\end{lem}
\begin{lem}\label{l2.12}Let $P_{*}(f)$ and $P_{*}(g)$ be defined in (\ref{e1.2}), for two non constant meromorphic functions $f$ and $g$. Then \beas &&\hspace{-15mm}\ol N(r,0;P_{*}(f)) \leq (\chi^{n-1}_{0}+m_{1}+m_{2}) T(r,f);\\&& \hspace{-15mm}N_{2}(r,0;P_{*}(f)) \leq (\chi^{n-1}_{0}+\mu^{n-1}_{0}+m_{1}+2m_{2}) T(r,f). \eeas
	Similar results occur for $P_{*}(g)$.
	\begin{proof}
	
 Rewrite $P_{*}(f)$ and $P_{*}(g)$ as \bea\label{e2.1} P_{*}(f)&=&C f(f-\beta_{_{1}}) \ldots(f-\beta_{_{m_{1}}})(f-\beta_{_{m_{1}+1}})^{n_{m_{1}+1}} \ldots (f-\beta_{_{m_{1}+m_{2}}})^{n_{m_{1}+m_{2}}}\eea and \beas P_{*}(g)&=&C g(g-\beta_{_{1}}) \ldots(g-\beta_{_{m_{1}}})(g-\beta_{_{m_{1}+1}})^{n_{m_{1}+1}} \ldots (g-\beta_{_{m_{1}+m_{2}}})^{n_{m_{1}+m_{2}}},\eeas where $\beta_{_{i}}'s$ $(i=1, 2, \ldots, m_{1}+m_{2})$ are distinct complex constants and $n_{i}$ is the multiplicity of the factor $(z-\beta_{_{i}})$ in $P^{*}(z)$ for $i=1, 2, \ldots, m_{1}+m_{2}$  with $n_{1}=n_{2}=\ldots=n_{m_{1}}=1$ and $n_{m_{1}+1}, \ldots, n_{m_{1}+m_{2}} \geq 2.$ \par Here  we have to consider two cases:\\
{\bf Case 1.} Suppose none of $\beta_{_{i}}'s$ $(i=1, 2, \ldots,m_{1}+m_{2})$ be zero. Then \beas \ol N(r,0;P_{*}(f)) \hspace{1mm}\leq\hspace{1mm}\ol N(r,0;f)+\displaystyle\sum_{i=1}^{m_{1}+m_{2}}\ol N(r,\beta_{_{i}};f) &\leq& (1+m_{1}+m_{2}) T(r,f) ;\eeas 
\beas  N_{2}(r,0;P_{*}(f)) &\leq& N(r,0;f)+\displaystyle\sum_{i=1}^{m_{1}} N(r,\beta_{_{i}};f)+ 2\hspace{-2mm}\displaystyle\sum_{i=m_{1}+1}^{m_{1}+m_{2}}\ol N(r,\beta_{_{i}};f) \hspace{.5cc}\leq (1+m_{1}+2m_{2}) T(r,f).\eeas
\textbf {Case 2:} Next let one of $\beta_{_{i}}'s$ $(i=1, 2, \ldots,m_{1}+m_{2})$ be zero.\\
Subcase 1: Suppose one among $\beta_{_{i}}'s$ $(i=1, 2, \ldots,m_{1})$ be zero. Without loss of generality let us assume that $\beta_{_{1}}=0$. Then
\beas \ol N(r,0;P_{*}(f)) &\leq& \ol N(r,0;f)+\displaystyle\sum_{i=2}^{m_{1}+m_{2}}\ol N(r,\beta_{_{i}};f) \hspace{1mm}\leq\hspace{1mm} (m_{1}+m_{2}) T(r,f) \;;\eeas 
\beas  N_{2}(r,0;P_{*}(f)) \hspace{.5mm}\leq \hspace{.5mm} 2 \ol N(r,0;f)+\displaystyle\sum_{i=2}^{m_{1}} N(r,\beta_{_{i}};f)+2 \hspace{-2mm}\displaystyle\sum_{i=m_{1}+1}^{m_{1}+m_{2}}\ol N(r,\beta_{_{i}};f) \hspace{.5mm}\leq \hspace{.5mm} (1+m_{1}+2m_{2}) T(r,f).\eeas 
Subcase 2: Next suppose one among $\beta_{_{i}}'s$ $(i=m_{1}+1, m_{1}+2, \ldots,m_{1}+m_{2})$ be zero. Without loss of generality let us assume that $\beta_{_{m_{1}+1}}=0$. Then
\beas \ol N(r,0;P_{*}(f)) &\leq& \ol N(r,0;f)+\displaystyle\sum_{i=1}^{m_{1}}\ol N(r,\beta_{_{i}};f)+\displaystyle\sum_{i=m_{1}+2}^{m_{1}+m_{2}}\ol N(r,\beta_{_{i}};f) \hspace{.5mm}\leq\hspace{.5mm} (m_{1}+m_{2}) T(r,f) ;\eeas 
\beas  N_{2}(r,0;P_{*}(f)) &\leq& 2 \ol N(r,0;f)+\displaystyle\sum_{i=1}^{m_{1}} N(r,\beta_{_{i}};f)+2\hspace{-2mm}\displaystyle\sum_{i=m_{1}+2}^{m_{1}+m_{2}}\ol N(r,\beta_{_{i}};f) \hspace{.5mm}\leq\hspace{.5mm} (m_{1}+2m_{2}) T(r,f).\eeas 
Combining all cases we can write 
\beas &&\ol N(r,0;P_{*}(f)) \leq (\chi^{n-1}_{0}+m_{1}+m_{2}) T(r,f);\\&& N_{2}(r,0;P_{*}(f)) \leq (\chi^{n-1}_{0}+\mu^{n-1}_{0}+m_{1}+2m_{2}) T(r,f). \eeas
Similarly we can obtain the same conclusions for the function $g$.	\end{proof}
\end{lem}
\begin{lem}\label{l2.13} Let $P_{*}(f)$ and $P_{*}(g)$ for two non constant meromorphic functions $f$ and $g$ (as defined in (\ref{e1.2})) share $(1,2)$ and $(\infty,0)$. If \beas n > 2(\chi^{n-1}_{0}+\mu^{n-1}_{0}+m_{1}+2m_{2})+\frac{15}{(2n-3)} (\chi^{n-1}_{0}+m_{1}+m_{2}), \eeas then either $P_{*}(f)(z)\equiv P_{*}(g)(z)$ or $P_{*}(f)(z).P_{*}(g)(z)\equiv 1$.
\begin{proof}
Set \beas \Phi=\frac{P_{*}(f)(P_{*}(g)-1)}{P_{*}(g)(P_{*}(f)-1)}.\eeas Clearly  $S(r,\Phi)$ can be replaced by $S(r,f)+S(r,g)$. It is obvious that $\Phi\not\equiv 0$. If $\Phi\equiv 0$ then either $P_{*}(f)=0$ or $P_{*}(g)=1$, which gives $f$ and $g$ are constants, a contradiction.\\
First suppose that $\Phi \not\equiv 1$. So 
$P_{*}(f)\not\equiv P_{*}(g)$.\\
Therefore, using {\it Lemma \ref{l2.10}} we get \beas&& \ol N(r,0;\Phi)+\ol N(r,\infty;\Phi)\\ &\leq& \ol N(r,1;P_{*}(f)\hspace{0.2cc} \mid \geq 3)+ \ol N(r,0;P_{*}(f))+\ol N(r,0;P_{*}(g)) \\&\leq& \frac{1}{2}\left(\ol N(r,0;P_{*}(f))+\ol N(r,\infty;P_{*}(f))\right) + \ol N(r,0;P_{*}(f))\\&&+\ol N(r,0;P_{*}(g))+S(r,P_{*}(f)) \\&\leq& \frac{3}{2}\ol N(r,0;P_{*}(f))+\frac{1}{2}\ol N(r,\infty;f)+\ol N(r,0;P_{*}(g))+S(r,f).\eeas
Now, \beas \Phi-1=\frac{P_{*}(g)-P_{*}(f)}{P_{*}(g)(P_{*}(f)-1)}\;\;\text{and}\;\; \Phi^{\prime}=\left[ \frac{P_{*}(g)^{\prime}}{P_{*}(g)(P_{*}(g)-1)}-\frac{P_{*}(f)^{\prime}}{P_{*}(f)(P_{*}(f)-1)}\right]\Phi. \eeas If $\Phi^{\prime}\equiv 0$ then \beas \left[ \frac{P_{*}(g)^{\prime}}{P_{*}(g)(P_{*}(g)-1)}-\frac{P_{*}(f)^{\prime}}{P_{*}(f)(P_{*}(f)-1)}\right]\equiv 0. \eeas Integrating we have, \beas \frac{P_{*}(f)-1}{P_{*}(f)}\equiv A \;\frac{P_{*}(g)-1}{P_{*}(g)}, \eeas where $A$ is non zero constant. i.e., \beas 1-\frac{1}{P_{*}(f)}\equiv A- \frac{A}{P_{*}(g)}.\eeas Since $P_{*}(f)$ and $P_{*}(g)$ share $(\infty, 0)$ so $A=1$. Then $P_{*}(f)\equiv P_{*}(g)$ which gives $\Phi\equiv 1,$ a contradiction. Therefore $\Phi^{}{\prime}\not\equiv 0$. Clearly all poles of $P_{*}(f)$ and $P_{*}(g)$ are multiple poles which are multiple  zeros of $\Phi-1$ and so zeros of $\Phi^{\prime}$ with multiplicity at least $(n-1)$  but not zeros of $\Phi$. Therefore by {\it Lemma \ref{l2.9}}, \beas(n-1) \ol N(r,\infty;f)&=&(n-1) \ol N(r,\infty;P_{*}(f)) = (n-1) \ol N(r,\infty;P_{*}(f)\mid \geq n)\nonumber \\&\leq& N(r,0;\Phi^{'} \mid \Phi \neq 0) \leq \ol N(r,0;\Phi)+\ol N(r,\infty;\Phi)+S(r,\Phi).\eeas
So, \beas (2n-3)\ol N(r,\infty;f)\leq 3\ol N(r,0;P_{*}(f))+2\ol N(r,0;P_{*}(g))+S(r,f).\eeas
 Applying {\it Lemma \ref{l2.12}} we obtain
\beas	\ol N(r,\infty;f)\leq \frac{3(\chi^{n-1}_{0}+m_{1}+m_{2})}{2n-3}T(r,f)+\frac{2(\chi^{n-1}_{0}+m_{1}+m_{2})}{2n-3} T(r,g)+S(r,f)+S(r,g).\eeas Similarly
\beas 	\ol N(r,\infty;g)\leq \frac{3(\chi^{n-1}_{0}+m_{1}+m_{2})}{2n-3}T(r,g)+\frac{2(\chi^{n-1}_{0}+m_{1}+m_{2})}{2n-3} T(r,f)+S(r,f)+S(r,g).\eeas That is \bea\label{e2.2}\ol N(r,\infty;f) + 	\ol N(r,\infty;g)&\leq& \frac{5(\chi^{n-1}_{0}+m_{1}+m_{2})}{2n-3}\left(T(r,f)+T(r,g)\right)\\&&+S(r,f)+S(r,g)\nonumber.\eea

If possible, we suppose that (i) of {\it Lemma \ref {l2.11}} holds. Therefore
\beas&& T(r,P_{*}(f))+T(r,P_{*}(g))\\&\leq& 2\{N_{2}(r,0;P_{*}(f))+N_{2}(r,0;P_{*}(g))+\ol N(r,\infty;P_{*}(f))+\ol N(r,\infty;P_{*}(g))\\&&+\ol N_{*}(r,\infty;P_{*}(f),P_{*}(g))\}+S(r,P_{*}(f))+S(r,P_{*}(g)).\eeas
Then using {\it Lemma \ref {l2.7}}, {\it Lemma \ref {l2.12}} and (\ref{e2.2}) we have \beas && n\left(T(r,f)+T(r,g)\right)\\&\leq& \left( 2(\chi^{n-1}_{0}+\mu^{n-1}_{0}+m_{1}+2m_{2})+\frac{15(\chi^{n-1}_{0}+m_{1}+m_{2})}{2n-3}\right) \left(T(r,f)+T(r,g)\right)\\&&+S(r,f)+S(r,g),\eeas which contradicts our assumption. So by {\it Lemma \ref {l2.11}} we have $$P_{*}(f)(z).P_{*}(g)(z)\equiv 1.$$ 
If $\Phi\equiv 1$, then $P_{*}(f)(z)\equiv P_{*}(g)(z).$
\par Hence the lemma is proved.
	\end{proof}
\end{lem}
\begin{lem}\label{l2.14} 
	Let $f$ and $g$ be two non constant meromorphic functions of finite order. Let $n\geq 2$, and let $ \{a_{1}(z)$, $a_{2}(z), \ldots, a_{n}(z)\} \in S(f)$ be distinct meromorphic periodic functions with period $c$. If $m\left(r, \frac{g}{f-a_{k}}\right)=S(r,f)$, for $k=1,2, \ldots,n$, then \beas \displaystyle\sum_{k=1}^{n} m\left(r, \frac{1}{f-a_{k}}\right)\leq m\left(r, \frac{1}{ g}\right)+S(r,f),\eeas where the exceptional set associated with $S(r,f)$ is of at most finite logarithmic measure.
	\begin{proof}
		Set \beas P(f)=\prod_{k=1}^{n}(f-a_{k}).\eeas Rewriting we have, \beas \frac{1}{P(f)}=\displaystyle\sum_{k=1}^{n}\frac{\alpha_{k}}{f-a_{k}}, \eeas where $\alpha_{k}\in S(f)$ are certain periodic function with period $c$. Now,\beas m\left (r,\frac{g}{P(f)}\right) \leq \displaystyle \sum_{k=1}^{n} m\left(r, \frac{g}{f-a_{k}}\right)+S(r,f)=S(r,f),\eeas  and so \beas m\left(r, \frac{1}{P(f)}\right)=m\left(r,\frac{g}{P(f)}\right)+ m\left(r, \frac{1}{g}\right)\leq m\left(r, \frac{1}{g}\right)+ S(r,f).\eeas By the first fundamental theorem and using the above inequation we get, 
		\beas m\left(r, \frac{1}{ g}\right)&\geq& m\left(r, \frac{1}{P(f)}\right)+
		S(r,f)= T(r, P(f))-N\left(r, \frac{1}{P(f)}\right)+S(r,f)\\&\geq& n T(r,f)-\displaystyle\sum_{k=1}^{n} N\left(r, \frac{1}{f-a_{k}}\right)+S(r,f)=\displaystyle\sum_{k=1}^{n} m\left(r, \frac{1}{f-a_{k}}\right)+S(r,f).\eeas 
	\end{proof} 
\end{lem}
\begin{lem}\label{l2.15} If $f$ be a meromorphic function of finite order then $\Tilde{L}(f(z))$ is of finite order and $m\left(r, \frac{\Tilde{L}(f(z))}{f(z+c)}\right)=S(r,f)$, $m\left(r, \frac{\Tilde{L}(f(z))}{f(z)-\beta_{i}}\right)=S(r,f)$ and $m\left(r, \frac{\Tilde{L}(f(z))}{f(qz)}\right)=S_{1}(r,f).$
	\begin{proof}
		 
		Using logarithmic derivative lemma and {\it Lemma \ref{l2.2}} we have,
		\bea\label{e2.4}\;\;\;\;\;\; m\left(r, \frac{\Tilde{L}(f(z))}{f(z+c)}\right)&=&m\left(r,\;\; \frac{\displaystyle\sum_{j=1}^{s}b_{j}f^{(j)}(z+c_{j})+\displaystyle\sum_{j=1}^{t}d_{j}f^{(j)}(z)}{f(z+c)}\right)\\&\leq&\displaystyle\sum_{j=1}^{s} m\left(r,\; \frac{f^{(j)}(z+c_{j})}{f^{(j)}(z)}\right)+\displaystyle\sum_{j=1}^{s} m\left(r, \frac{f^{(j)}(z)}{f(z)}\right)\nonumber\\&&+\displaystyle\sum_{j=1}^{t} m\left(r, \frac{f^{(j)}(z)}{f(z)}\right)+(s+t)\;m\left(r, \frac{f(z)}{f(z+c)}\right)+O(1)\nonumber\\&=&S(r,f)\nonumber.\eea
		Also, \beas m\left(r, \frac{\Tilde{L}(f(z))}{f(z)-\beta_{i}}\right)&=&m\left(r,\; \frac{\displaystyle\sum_{j=1}^{s}b_{j}f^{(j)}(z+c_{j})+\displaystyle\sum_{j=1}^{t}d_{j}f^{(j)}(z)}{f(z)-\beta_{i}}\right)\\&\leq&\displaystyle\sum_{j=1}^{s} m\left(r, \frac{f^{(j)}(z+c_{j})}{f^{(j)}(z)}\right)+\displaystyle\sum_{j=1}^{t} m\left(r, \frac{f^{(j)}(z)}{f(z)-\beta_{i}}\right)\nonumber\\&&+\displaystyle\sum_{j=1}^{s} m\left(r, \frac{f^{(j)}(z)}{f(z)-\beta_{i}}\right)+O(1)=S(r,f)\nonumber.\eeas
		Using $(\ref{e2.4})$ and {\it Lemma \ref{l2.1}} we have, \beas T(r, \Tilde{L}(f(z)))\leq \frac{s^{2}+t^{2}+3(s+t)+2}{2}\;T(r, f)+S(r,f) . \eeas As $f$ is  of finite order so $\Tilde{L}(f(z))$ and $f(z+c)$ is of finite order and $S(r, \Tilde{L}(f(z)))$ can be replaced by $S(r,f)$.\par Similarly by using {\it Lemma \ref{l2.4}}, {\it Lemma \ref{l2.5}}  and {\it Lemma \ref{l2.6}} as and when required we can prove  $f(qz)$ and $\Tilde{L}(f(z))$ are zero order when $f$ is of zero order and
		\beas m\left(r, \frac{\Tilde{L}(f(z))}{f(qz)}\right)=S_{1}(r,f).\eeas
		
	\end{proof} 
\end{lem}
\section {Proofs of the theorems} 
\begin{proof}[ {\bf Proof of Theorem \ref{t1}}] 
Since $E_{f(z+c)}(S^{*}_{1},2)=E_{\Tilde{L}(f(z))}(S^{*}_{1},2)$ and $E_{f(z+c)}(S_{2},0)=E_{\Tilde{L}(f(z))}\\(S_{2},0)$, it follows that $P_{*}(f(z+c))$, $P_{*}(\Tilde{L}(f(z)))$ share $(1,2)$ and $(\infty,0)$.
So by {\it Lemma \ref {l2.13}} we have either $P_{*}(f(z+c))\equiv P_{*}(\Tilde{L}(f(z)))$ or $P_{*}(f(z+c)).P_{*}(\Tilde{L}(f(z)))\equiv 1.$ Suppose that \bea\label{e3.1}P_{*}(f(z+c)).P_{*}(\Tilde{L}(f(z)))\equiv 1.\eea Noting that $P_{*}(f(z+c))$, $P_{*}(\Tilde{L}(f(z)))$ share $(\infty,0)$, so we can conclude that $P_{*}(f(z+c))$, $P_{*}(\Tilde{L}(f(z)))$ both are entire functions.\par So $$N\left(r,\infty;\frac{P_{*}(\Tilde{L}(f(z)))}{P_{*}(f(z+c))}\right)= N(r,0;P_{*}(f(z+c))).$$\par 
Therefore using {\it Lemma \ref {l2.12}} and {\it Lemma \ref {l2.1}}, we get \beas N\left(r,\infty;\frac{P_{*}(\Tilde{L}(f(z)))}{P_{*}(f(z+c))}\right)\leq (\chi^{n-1}_{0}+m_{1}+m_{2}) T(r,f(z+c))\leq nT(r,f)+S(r,f).\eeas 
Using {\it Lemma \ref{l2.2}} and {\it Lemma \ref{l2.15}} we have, \beas m\left(r,\;\frac{P_{*}(\Tilde{L}(f(z)))}{P_{*}(f(z+c))}\right)&=&m\left(r,\; \frac{\Tilde{L}(f(z))}{f(z+c)}\prod_{i=1}^{m_{1}+m_{2}}\left(\frac{\Tilde{L}(f(z))-\beta_{_{i}}}{f(z+c)-\beta_{_{i}}}\right)^{n_{i}}\right)\\&\leq& m\left(r, \frac{\Tilde{L}(f(z))}{f(z+c)}\right)+ m\left(r,\prod_{i=1
}^{m_{1}+m_{2}}\left(\frac{\Tilde{L}(f(z))-\beta_{_{i}}}{f(z+c)-\beta_{_{i}}}\right)^{n_{i}}\right)+O(1)\\&\leq& \displaystyle
\sum_{i=1}^{m_{1}+m_{2}} n_{i} \;m\left(r,\frac{\Tilde{L}(f(z))-\beta_{_{i}}}{f(z+c)-\beta_{_{i}}}\right)+S(r,f)\\&\leq& \displaystyle
\sum_{i=1}^{m_{1}+m_{2}} n_{i} \;m\left(r,\frac{\Tilde{L}(f(z))}{f(z)-\beta_{_{i}}}\right)+\displaystyle
\sum_{i=1}^{m_{1}+m_{2}} n_{i} \;m\left(r,\frac{1}{f(z)-\beta_{_{i}}}\right)\\&&+\displaystyle
\sum_{i=1}^{m_{1}+m_{2}} n_{i} \;m\left(r,\frac{f(z)-\beta_{_{i}}}{f(z+c)-\beta_{_{i}}}\right)+S(r,f)\\&\leq&\displaystyle
\sum_{i=1}^{m_{1}+m_{2}} n_{i} \;m\left(r,\frac{1}{f(z)-\beta_{_{i}}}\right)+S(r,f)\\&\leq&(n_{1}+n_{2}+\ldots+n_{m_{1}+m_{2}})T(r,f)+S(r,f)\\&\leq&(n-1)T(r,f)+S(r,f).\eeas By {\it Lemma \ref{l2.1}}, {\it Lemma \ref{l2.7}} and $(\ref{e3.1})$,
\beas 2nT(r,f)&=& 2nT(r,f(z+c))+S(r,f)=2T(r,P_{*}(f(z+c)))+S(r,f)\\&\leq& T\left(r,\;\frac{1}{P_{*}(f(z+c))^{2}}\right)+S(r,f)\leq T\left(r,\;\frac{P_{*}(\Tilde{L}(f(z)))}{P_{*}(f(z+c))}\right)+S(r,f)\\&\leq& (2n-1)T(r,f)+S(r,f),\eeas which is a contradiction.\par
Therefore $P_{*}(\Tilde{L}(f(z)))\equiv P_{*}(f(z+c))$, which yields \beas \prod_{i=1}^{n} (\Tilde{L}(f(z))-\alpha_{i}) \equiv \prod_{i=1}^{n} (f(z+c)-\alpha_{i}).\eeas\end{proof} 
\begin{proof}[ {\bf Proof of Theorem \ref{t2}}] 
By proceeding in a similar way of the proof of {\it Theorem \ref{t1}} we can prove this theorem using {\it Lemma \ref {l2.4}}, {\it Lemma \ref {l2.5}} and {\it Lemma \ref {l2.6}} as and when required instead of {\it Lemma \ref {l2.1} and }{\it Lemma \ref {l2.2}}.
\end{proof} 
\begin{proof}[ {\bf Proof of Theorem \ref{t3}}]
 Since the finite order meromorphic functions $f(z+c)$ and $\Tilde{L}(f(z))$ share $(S^{*}_{1}, \infty)$, $(S_{2}, \infty)$, it follows that  $P_{*}(f(z+c))$, $P_{*}(\Tilde{L}(f(z)))$ share $(1, \infty)$ and $(\infty,\infty)$ which yields \bea\label{e3.2} N(r,\Tilde{L}(f(z)))=N(r,f(z+c)) \eea and \bea\label{e3.3} \frac{P_{*}(\Tilde{L}(f(z)))-1}{P_{*}(f(z+c))-1}=e^{\gamma(z)},\eea where $\gamma(z)$ is a polynomial.\par Now, \beas T(r, e^{\gamma(z)})=m(r,e^{\gamma(z)})=m\left(r,\frac{P_{*}(\Tilde{L}(f(z)))-1}{P_{*}(f(z+c))-1}\right).\eeas Using the definition of $P_{*}(z)$ we have,
	\beas T(r, e^{\gamma(z)})&=&m\left(r,\frac{(\Tilde{L}(f(z)))-\alpha_{1})(\Tilde{L}(f(z)))-\alpha_{2})\ldots(\Tilde{L}(f(z)))-\alpha_{n})}{(f(z+c)-\alpha_{1})(f(z+c)-\alpha_{2})\ldots(f(z+c)-\alpha_{n})}\right)\\&\leq& \displaystyle\sum_{j=1}^{n}m\left(r,\frac{\Tilde{L}(f(z)))-\alpha_{j}}{f(z+c)-\alpha_{j}}\right)+O(1)\\&\leq&\displaystyle\sum_{j=1}^{n}m\left(r,\frac{\Tilde{L}(f(z))}{f(z)-\alpha_{j}}\right)+\displaystyle\sum_{j=1}^{n}m\left(r,\frac{1}{f(z)-\alpha_{j}}\right)+\displaystyle\sum_{j=1}^{n}m\left(r,\frac{f(z)-\alpha_{j}}{f(z+c)-\alpha_{j}}\right)\\&&+O(1) .\eeas In view of {\it Lemma \ref{l2.2}}, {\it Lemma \ref{l2.14}}, {\it Lemma \ref{l2.15}} and then by the first fundamental theorem and $(\ref{e3.2})$ we have, \beas T(r, e^{\gamma(z)})&=&\displaystyle\sum_{j=1}^{n}m\left(r,\frac{1}{f(z)-\alpha_{j}}\right)+S(r,f)\leq m\left(r,\frac{1}{\Tilde{L}(f(z))}\right)+S(r,f)\\&\leq& T(r,\Tilde{L}(f(z)))-N\left(r,\frac{1}{\Tilde{L}(f(z))}\right)+S(r,f)\\&\leq& m\left(r, \frac{\Tilde{L}(f(z))}{f(z+c)}\right)+m(r,f(z+c))+N(r,\Tilde{L}(f(z)))-N\left(r,\frac{1}{\Tilde{L}(f(z))}\right)+S(r,f)\\&\leq& T(r,f(z+c))-N\left(r,\frac{1}{\Tilde{L}(f(z))}\right)+S(r,f)\leq T(r,f)-N\left(r,\frac{1}{\Tilde{L}(f(z))}\right)+S(r,f).\eeas
	According to the given condition $T(r,f)=N\left(r,\frac{1}{ \Tilde{L}(f(z))}\right)+S(r,f)$, so $T(r, e^{\gamma(z)})=S(r,f).$\par
	Now from $(\ref{e3.3})$ we have, \beas P_{*}(\Tilde{L}(f(z)))= e^{\gamma(z)}\left(P_{*}(f(z+c))-1+e^{-\gamma(z)}\right). \eeas	Set \beas W(z)=\frac{P_{*}(f(z+c))}{1-e^{-\gamma(z)}} .\eeas	If $e^{\gamma(z)}\not\equiv 1$, then by applying Nevanlinna's second fundamental theorem to $W(z)$ and using $(\ref{e3.2})$ and {\it Lemma \ref{l2.12}} we obtain,	\beas& & T(r,P_{*}(f(z+c)))\leq T(r,W)+S(r,f)\\&\leq& \ol N(r,0; W)+ \ol N(r,\infty;W)+ \ol N(r,0; W-1)+S(r,f)\\&\leq& \ol N(r,0;P_{*}(f(z+c)))+ \ol N(r,\infty;P_{*}(f(z+c))) + \ol N(r,0; P_{*}(\Tilde{L}(f(z))))+ S(r,f)\\&\leq& (\chi^{n-1}_{0}+m_{1}+m_{2})\left(T(r,f(z+c))+T(r,\Tilde{L}(f(z)))\right)+ N(r,\infty;f(z+c))+S(r,f) \\&\leq& (\chi^{n-1}_{0}+m_{1}+m_{2})\left(T(r,f(z+c))+m(r,f(z+c))+m\left(r,\frac{\Tilde{L}(f(z))}{f(z+c)}\right)\right.\\&&\left.+N(r,\infty;f(z+c))\right)+N(r,\infty;f)+S(r,f).\eeas	Using {\it Lemma \ref{l2.1}} and {\it Lemma \ref{l2.15}} we get, \beas nT(r,f)\leq (2\chi^{n-1}_{0}+2m_{1}+2m_{2}+1)T(r,f)+S(r,f),\eeas which contradicts $n > 2(\chi^{n-1}_{0}+m_{1}+m_{2})+1$. This gives $e^{\gamma(z)}\equiv 1$, that yields \beas \prod_{i=1}^{n} (\Tilde{L}(f(z))-\alpha_{i}) \equiv \prod_{i=1}^{n} (f(z+c)-\alpha_{i}).\eeas
\end{proof} 
\begin{proof}[ {\bf Proof of Theorem \ref{t4}}] Here $\Tilde{L}(f(z))$ and $f(qz)$ are of zero order. Since $f(qz)$ and $\Tilde{L}(f(z))$ share $(S^{*}_{1}, \infty)$ and $(S_{2}, \infty)$, it follows that  $P_{*}(f(qz))$ and $P_{*}(\Tilde{L}(f(z)))$ share $(1, \infty)$ and $(\infty,\infty)$. Therefore \beas \frac{P_{*}(\Tilde{L}(f(z)))-1}{P_{*}(f(qz))-1}=A,\eeas where A is a non zero constant.\par	This gives \beas P_{*}(\Tilde{L}(f(z)))= A\left(P_{*}(f(qz))-1+\frac{1}{A}\right). \eeas	Set $ W_{1}(z)=\frac{P_{*}(f(qz))}{1-\frac{1}{A}}. $	If $A\not\equiv 1$, then applying Nevanlinna's second fundamental theorem to $W_{1}(z)$ and using  {\it Lemmas \ref{l2.4}} and {\it \ref{l2.5}} and {\it \ref{l2.15}} as and when required we can calculate the rest of the proof similar to {\it Theorem \ref{t3}}.	 \end{proof}

\begin{proof}[ {\bf Proof of Theorem \ref{t5}}]Here $f(z+c)$ and $\Tilde{L}(f(z))$ are of finite order. Since $f(z+c)$ and $\Tilde{L}(f(z))$ share $(0, \infty)$ and $(\infty, \infty)$, so \bea\label{e3.4} \frac{\Tilde{L}(f(z))}{f(z+c)}=e^{\delta(z),}\eea where $\delta(z)$	is a polynomial. \\Clearly by {\it Lemma \ref{l2.15}} we get, \beas T(r,e^{\delta(z)})=S(r,f).\eeas When $e^{\delta(z)}\equiv 1$ then $\Tilde{L}(f(z))\equiv f(z+c)$. \par When $e^{\delta(z)}\not\equiv 1$, using the fact that $f(z+c)$ and $\Tilde{L}(f(z))$ share $(a(z), 0)$ we have, \beas \ol N\left(r, \frac{1}{\Tilde{L}(f(z))-a(z)}\right)&=&\ol N\left(r, \frac{1}{f(z+c)-a(z)}\right)\leq\ol N\left(r, \frac{1}{e^{\delta(z)}-1}\right)+\ol N\left(r, \frac{1}{a(z)}\right)\\&\leq& T(r, e^{\delta(z)})+S(r,f)=S(r,f).\eeas Rewriting $(\ref{e3.4})$ we get, \beas \Tilde{L}(f(z))-a(z)=e^{\delta(z)}(f(z+c)-a(z)e^{-\delta(z)}).  \eeas Clearly $a(z)e^{-\delta(z)}\not\equiv a(z)$. So, \beas\ol N\left(r, \frac{1}{f(z+c)-a(z)e^{-\delta(z)}}\right)=\ol N\left(r, \frac{1}{\Tilde{L}(f(z))-a(z)}\right)=S(r,f).\eeas Using {\it Lemma \ref{l2.1}}, {\it \ref{l2.3}} and the second fundamental theorem we obtain, \beas&& 2T(r,f)=2T(r, f(z+c))+S(r,f)\\&\leq& \ol N(r, f(z+c))+\ol N\left(r, \frac{1}{f(z+c)}\right)+\ol N\left(r, \frac{1}{f(z+c)-a(z)}\right)\\&&+\ol N\left(r, \frac{1}{f(z+c)-a(z)e^{-\delta(z)}}\right)+S(r,f) \leq\ol N(r, f)+\ol N\left(r, \frac{1}{f}\right)+S(r,f),\eeas which is a contradiction to $\Theta(0;f)+\Theta(\infty; f)> 0$. Hence $\Tilde{L}(f(z))\equiv f(z+c)$.
 \end{proof}
\begin{proof}[ {\bf Proof of Theorem \ref{t6}}] Here $f(z+c)$ and $\Tilde{L}(f(z))$ are of finite order. Since $f(z+c)$ and $\Tilde{L}(f(z))$ share $(a(z), \infty)$ and $(\infty, \infty)$, so \bea\label{e3.5} \frac{\Tilde{L}(f(z))-a(z)}{f(z+c)-a(z)}=e^{\zeta(z),}\eea where $\zeta(z)$	is a polynomial. Using logarithmic derivative lemma, {\it Lemma \ref{l2.1}} and {\it Lemma \ref{l2.2}} we get,   \beas &&T(r,e^{\zeta(z)})=m(r,e^{\zeta(z)})=m\left(r,\; \frac{\Tilde{L}(f(z))-a(z)}{f(z+c)-a(z)}\right)\\&\leq& m\left(r,\; \frac{\Tilde{L}(f(z))-\Tilde{L}(a(z-c))}{f(z+c)-a(z)}\right)+m\left(r,\; \frac{\Tilde{L}(a(z-c))-a(z)}{f(z+c)-a(z)}\right)\\&\leq&m\left(r,\; \frac{\Tilde{L}(f(z))-\Tilde{L}(a(z-c))}{f(z)-a(z-c)}\right)+m\left(r,\; \frac{f(z)-a(z-c)}{f(z+c)-a(z)}\right)\\&&+m\left(r,\; \frac{1}{f(z+c)-a(z)}\right)+S(r,f)\\&\leq&m\left(r,\; \frac{\displaystyle\sum_{j=1}^{s}b_{j}(f^{(j)}(z+c_{j})-a^{(j)}(z-c+c_{j}))+\displaystyle\sum_{j=1}^{t}d_{j}(f^{(j)}(z)-a^{(j)}(z-c))}{f(z)-a(z-c)}\right)\\&&+T(r, f(z+c))+S(r,f)\\&\leq&\displaystyle\sum_{j=1}^{s} m\left(r,\; \frac{f^{(j)}(z+c_{j})-a^{(j)}(z-c+c_{j})}{f^{(j)}(z)-a^{(j)}(z-c)}\right)+\displaystyle\sum_{j=1}^{t} m\left(r,\; \frac{f^{(j)}(z)-a^{(j)}(z-c)}{f(z)-a(z-c)}\right)\\&&+\displaystyle\sum_{j=1}^{s} m\left(r,\; \frac{f^{(j)}(z)-a^{(j)}(z-c)}{f(z)-a(z-c)}\right)+T(r,f)+S(r,f)\nonumber\\&\leq&T(r, f)+S(r,f).\eeas So $S(r,e^{\zeta(z)})$ can be replaced by $S(r,f)$. When $e^{\zeta(z)}\equiv 1$ then  $\Tilde{L}(f(z))\equiv f(z+c)$. \\Suppose $e^{\delta(z)}\not\equiv 1$. Now rewriting $(\ref{e3.5})$ we can obtain, \beas \frac{1}{f(z+c)}=-\frac{\Tilde{L}(f(z))}{f(z+c)(e^{\zeta(z)}-1)}+\frac{e^{\zeta(z)}}{a(z)(e^{\zeta(z)}-1)}. \eeas Therefore in view of {\it Lemma \ref{l2.15}} we have, \beas m\left(r, \frac{1}{f(z+c)}\right)\leq 2 \;m\left(r, \frac{1}{e^{\zeta(z)}-1}\right)+S(r,f).\eeas If $\zeta(z)$ is constant then automatically $m\left(r, \frac{1}{f(z+c)}\right)=S(r,f)$. If $\zeta(z)$ is non constant then by {\it Lemma \ref{l2.8}} we get, $$m\left(r, \frac{1}{f(z+c)}\right)=S(r,e^{\zeta(z)})=S(r,f).$$ By {\it Lemma \ref{l2.1}} and {\it Lemma \ref{l2.3}} we have,\beas T(r,f)&=&T(r, f(z+c))+S(r,f)=T\left(r, \frac{1}{f(z+c)}\right)+S(r,f)\\&\leq& N\left(r, \frac{1}{f(z+c)}\right)+S(r,f) \leq N\left(r, \frac{1}{f}\right)+S(r,f)\leq T(r,f)+S(r,f).\eeas  Therefore,  \beas N\left(r, \frac{1}{f}\right)=T(r,f)+S(r,f),\eeas which contradicts the fact that $\delta(0,f)> 0$. Hence $\Tilde{L}(f(z))\equiv f(z+c)$.
\end{proof}

\section {Observation}
Take $\Tilde{L}(f(z))=L_{3}$ with all coefficients are $1$. Then we see that choosing $c=\frac{\log(\alpha+\alpha^{2}+\ldots+\alpha^{t})}{\alpha}$, where $1+\alpha+\ldots+\alpha^{t-1}\neq 0$, we somehow get a solution $f(z)=e^{\alpha z}$ $(\alpha\neq0)$ of \bea\label{e4.1}\Tilde{L}(f(z))=f(z+c). \eea However choosing $c=\frac{\pi}{2}$, we can present the solution of $f^{\prime}=f(z+c)$ as the linear combination of two independent solutions. e.g., $f(z)=d_{1}e^{iz}+d_{2}e^{-iz}$. So it is a matter of concern that how the solutions of (\ref{e4.1}) looks like. Unfortunately we can not elucidated in this matter.

\end{document}